# Finite-sum Composition Optimization via Variance Reduced Gradient Descent

Xiangru Lian[*], Mengdi Wang[†], and Ji Liu[*]
[*]Department of Computer Science, University of Rochester
[†]Department of Operations Research and Financial Engineering, Princeton University
xiangru@yandex.com, mengdiw@princeton.edu, ji.liu.uwisc@gmail.com

May 20, 2017

**Abstract**

*The stochastic composition optimization proposed recently by Wang et al. [2014] minimizes the objective with the compositional expectation form: $\min_x \ (\mathbb{E}_i F_i \circ \mathbb{E}_j G_j)(x)$. It summarizes many important applications in machine learning, statistics, and finance. In this paper, we consider the finite-sum scenario for composition optimization:*

$$\min_x f(x) := \frac{1}{n} \sum_{i=1}^{n} F_i \left( \frac{1}{m} \sum_{j=1}^{m} G_j(x) \right).$$

*We propose two algorithms to solve this problem by combining the stochastic compositional gradient descent (SCGD) and the stochastic variance reduced gradient (SVRG) technique. A constant linear convergence rate is proved for strongly convex optimization, which substantially improves the sublinear rate $O(K^{-0.8})$ of the best known algorithm.*

## 1 Introduction

The stochastic composition optimization proposed recently by Wang et al. [2014] minimizes the objective with the compositional expectation form:

$$\min_x \ (\mathbb{E}_i F_i \circ \mathbb{E}_j G_j)(x).$$

It has many emerging applications, ranging from machine and reinforcement learning [Dai et al., 2016, Wang et al., 2016] to risk management [Dentcheva et al., 2015]. It is also related to multi-stage stochastic programming [Shapiro et al., 2014] and adaptive simulation [Hu et al., 2014].

In general the stochastic composition optimization is substantially more difficult than the traditional stochastic optimization: $\min_x \mathbb{E}_i F_i(x)$. This is because the composition objective is no longer linear with respect to the joint distribution of data indices $(i, j)$. For example, the best-known algorithms studied in [Wang et al., 2014, 2016] achieve a finite-sample error bound $O(K^{-0.8})$ for strongly convex composition optimization, which deteriorates from the optimal rate $O(K^{-1})$ for the generic stochastic optimization.

In this paper, we study the *finite-sum* scenario for composition optimization in the following form

$$\min_{x \in \mathbb{R}^N} f(x) := F \circ G(x) = F(G(x)), \tag{1}$$

where the inner function $G : \mathbb{R}^N \to \mathbb{R}^M$ is the empirical mean of $m$ component functions $G_i : \mathbb{R}^N \to \mathbb{R}^M$:

$$G(x) = \frac{1}{m} \sum_{i=1}^{m} G_i(x),$$

and the outer function $F : \mathbb{R}^M \to \mathbb{R}$ is the empirical mean of $n$ component functions $F_j : \mathbb{R}^M \to \mathbb{R}$:

$$F(y) = \frac{1}{n} \sum_{j=1}^{n} F_j(y).$$

The finite-sum composition problem models optimization involving two fixed-size empirical data sets. Randomly selecting component functions $G_i, F_j$ can be viewed as randomized retrieval from each of the two data sets.

In this paper, we propose two efficient algorithms (namely, compositional SVRG-1 and compositional SVRG-2) for the finite-sum composition optimization problem (1). The new algorithms are developed by combining the stochastic compositional gradient descent (SCGD) technique [Wang et al., 2014, 2016] and the stochastic variance reduced gradient (SVRG) technique [Johnson and Zhang, 2013]. The new algorithms are motivated by the fact that SVRG is able to improve the sublinear convergence rate of stochastic gradient descent to linear convergence in the case of classical (strongly convex) finite-sum stochastic optimization. We prove that the two algorithms converge linearly for the finite-sum stochastic composition optimization, with query complexity $O\left((m + n + \tilde{\kappa}_1^4) \log(1/\epsilon)\right)$ and $O\left((m + n + \tilde{\kappa}_2^3) \log(1/\epsilon)\right)$ respectively (the $\tilde{\kappa}_1$ and $\tilde{\kappa}_2$ are two variants of the condition number, and their definitions will be specified later). To the best of our knowledge, this is the first work on finite-sum stochastic composition optimization and linearly convergent algorithms.

## 1.1 Related Works

This section reviews algorithms related to composition optimization and SVRG.

*Composition Optimization* draws much attention recently. Contrary to classical stochastic problems, the objective in composition optimization is no longer a plain summation of component functions.

Table 1: Convergence rates of stochastic composition gradient descent and SVRG. $K$ is the number of iterations. For fair comparison, the convergence rates for the proposed algorithms (Compositional SVRG-1 and Compositional SVRG-2) have taken the query complexity per iteration (epoch) into consideration.

|  | Nonconvex | Convex | Strongly Convex |
|---|---|---|---|
| Basic SCGD [Wang et al., 2014] | $O(K^{-1/4})$ | $O(K^{-1/4})$ | $O(K^{-2/3})$ |
| Accelerating SCGD [Wang et al., 2014] | $O(K^{-2/7})$ | $O(K^{-2/7})$ | $O(K^{-4/5})$ |
| Accelerating SCGD [Wang et al., 2016] <br> ∗ means if $G_j(\cdot)$'s are linear | $O(K^{-4/9})$ <br> or <br> $O(K^{-1/2})^*$ | $O(K^{-2/7})$ <br> or <br> $O(K^{-1/2})^*$ | $O(K^{-4/5})$ <br> or <br> $O(K^{-1})^*$ |
| SVRG ($0 < \rho < 1$) | - | - | $O\left(\rho^{\frac{K}{m+n+m\kappa}}\right)$ |
| Compositional SVRG-1 ($0 < \rho < 1$) | - | - | $O\left(\rho^{\frac{K}{m+n+\bar{\kappa}_1^4}}\right)$ |
| Compositional SVRG-2 ($0 < \rho < 1$) | - | - | $O\left(\rho^{\frac{K}{m+n+\bar{\kappa}_2^3}}\right)$ |

Given an index $i \in \{1, 2, \ldots, n\}$, we cannot use a single query to the oracle to get the gradient of a component function $F_i$ with respect to the optimization variable $x$. In contrast, we can query the gradient of $F_i$ with respect to an intermediate variable $G$ in a single query, and $G$ is itself a summation of $m$ component functions. Thus to calculate the gradient of a single component function in classical stochastic algorithms, we need at least $O(m)$ queries. When $m$ becomes large, the query complexity for classical stochastic optimization algorithms will significantly increase. This encourages people to search for a more sophisticated way to solve such problems.

Wang et al. [2014] proposed the generic composition optimization for the first time. Two stochastic algorithms - Basic SCGD and accelerating SCGD - are proposed for such optimization, with provable convergence rates. A recent work by Wang et al. [2016] improves the convergence rate of accelerating compositional SGD and finds that the optimal convergence rate can be obtained if $G_j(\cdot)$'s are linear. All convergence rates together with traditional SVRG are listed in Table 1. Note that for sufficient large $m$, the traditional SVRG algorithm will be slower than Compositional SVRG algorithms. In addition, some special cases such as risk optimization are studied in [Dentcheva et al., 2015].

*SVRG* is a very powerful technique for large scale optimization. This variance reduced optimization algorithm was originally developed in Johnson and Zhang [2013]. Its main advantage lies on its low storage requirement compared with other variance reduced algorithms [Defazio et al., 2014a,b, Schmidt et al., 2013]. The SVRG technique has been extended by many works [Allen-Zhu and Yuan, 2015, Harikandeh et al., 2015, Kolte et al., 2015, Konečný et al., 2016, Nitanda, 2014, 2015, Xiao and Zhang, 2014] for solving stochastic optimization. Similar algorithms include Konecnỳ and Richtárik [2013], Shalev-Shwartz and Zhang [2013]. For SVRG applied on classical stochastic problems, Johnson and Zhang [2013] proved a $O(\rho^s)$ convergence rate on strongly convex objectives, where $\rho$ is a constant smaller than 1 and $s$ is the epoch number. The query complexity per epoch is $n + \kappa$ where $n$ is the number of component functions and $\kappa$ is the condition number of the objective. For classical gradient descent, to obtain the same convergence rate, the query complexity per iteration has to be $n\kappa$ (see Nesterov [2013]), which is generally larger than $n + \kappa$. We list the results for SVRG on classical stochastic optimization problems with various types of objectives in Table 2. Gong and Ye [2014] extends the analysis of SVRG for strongly convex optimization to the more general *optimally* strongly convex optimization [Liu and Wright, 2015, Wang et al., 2014] and proves similar convergence rate. Recently, the asynchronous parallel version of SVRG is also studied [Reddi et al., 2015].

Table 2: Query complexity and convergence rate for SVRG on different objectives. QCPE stands for query complexity per epoch. $\rho$ is a constant smaller than 1. $s$ is the epoch number. $n$ is the number of component functions.

| SVRG | Nonconvex [Reddi et al., 2016] [Allen-Zhu and Hazan, 2016] | Convex [Reddi et al., 2016] | Strongly Convex [Johnson and Zhang, 2013] |
|---|---|---|---|
| Rate | $O(1/(sn^{1/3}))$ | $O(1/(s\sqrt{n}))$ | $O(\rho^s)$ |
| QCPE | $O(n)$ | $O(n)$ | $O(n + \kappa)$ |

## 1.2 Notation

Throughout this paper, we use the following simple notations.

- $[z]_i$ denotes the $i$th component of vector (or vector function) $z$;
- We use the following notations for derivatives of functions. Given any smooth function

$$H : \mathbb{R}^N \to \mathbb{R}^M$$
$$x \mapsto H(x),$$

$\partial H$ is the Jacobian of $H$ defined by

$$\partial H := \frac{\partial H}{\partial x} = \begin{pmatrix} \frac{\partial [H]_1}{\partial [x]_1} & \cdots & \frac{\partial [H]_1}{\partial [x]_N} \\ \vdots & \ddots & \vdots \\ \frac{\partial [H]_M}{\partial [x]_1} & \cdots & \frac{\partial [H]_M}{\partial [x]_N} \end{pmatrix}.$$

The value of the Jacobian at some point $a$ is denoted by $\partial H(a)$.

For a scalar function

$$h : \mathbb{R}^N \to \mathbb{R}$$
$$x \mapsto h(x),$$

the gradient of $h$ is defined by

$$\nabla h(x) = \left(\frac{\partial h(x)}{\partial x}\right)^\top = \left(\frac{\partial h}{\partial [x]_1}, \ldots, \frac{\partial h}{\partial [x]_N}\right)^\top \in \mathbb{R}^N.$$

Under this set of notations, with the chain rule we can easily find the gradient of a composition function $f = F(G(x))$ to be:

$$\nabla f(x) = (\partial G(x))^\top \nabla F(G(x)).^1$$

- $x^*$ denotes the optimal solution of (1);

- Given a multiset $\mathcal{A}^2$, we use $\text{len}(\mathcal{A})$ to denote the number of elements in $\mathcal{A}$. For example, if $\mathcal{A} = \{1, 2, 3, 1\}$, then $\text{len}(\mathcal{A}) = 4$. We use $\mathcal{A}[i]$ to represent the $i$th element in $\mathcal{A}$.

- $\mathbb{E}_i$ denotes taking expectation w.r.t. the random variable $i$.

- $\mathbb{E}$ denotes taking expectation w.r.t. all random variables.

## 2 Preliminary: SVRG

We review the standard SVRG algorithm in this section for completion. Consider to solve the following finite sum optimization problem

$$\min_x f(x) = \frac{1}{n} \sum_{i=1}^{n} F_i(x).$$

The SVRG algorithm basically stores the gradient of $f$ at a reference point $\tilde{x}$ (the reference point will be updated for every a few iterations): $\tilde{f}' := \frac{1}{n} \sum_{i=1}^{n} \nabla F_i(\tilde{x})$. Based on such a reference gradient, SVRG estimates the gradient at each iteration by

$$\hat{f}'_k := \tilde{f}' - \nabla F_i(\tilde{x}) + \nabla F_i(x_k)$$

where $i$ is uniformly randomly sampled from $\{1, 2, \ldots, n\}$. The next iterate is updated by

$$x_{k+1} = x_k - \gamma_k \hat{f}'_k.$$

The computation complexity per iteration is comparable to SGD. The estimated gradient is also an *unbiased* estimate for the true gradient

$$\mathbb{E}(\hat{f}'_k) = f'_k := \frac{1}{n} \sum_{i=1}^{n} \nabla F_i(x_k).$$

The key improvement lies on that the variance $\mathbb{E}(\|\hat{f}'_k - f'_k\|^2)$ decreases to zero when $x_k$ converges to the optimal point for SVRG while it is a constant for SGD. Therefore, SVRG admits a much better convergence rate (linear convergence for strongly convex optimization and sublinear for convex optimization) than SGD. For completeness, the complete SVRG algorithm is shown in Algorithm 1.

---

[1]Note that the gradient operator always calculates the gradient with respect to the first level variable. That is to say, by $\nabla F(G(x))$ we mean the gradient of $F(y)$ at $y = G(x)$, *not* the gradient of $F(G(x))$ with respect to $x$.

[2]A multiset is a generalization of the concept of a set, allowing duplicated elements

**Algorithm 1** SVRG [Johnson and Zhang, 2013]
---
**Require:** $K$ (update frequency), $\gamma$ (step length), $S$ (total number of epochs), $\tilde{x}_0$ (initial point)
**Ensure:** $\tilde{x}_S$.
 1: **for** $s = 1, 2, \ldots, S$ **do**
 2:     Update the reference point $\tilde{x} \leftarrow \tilde{x}_{s-1}$
 3:     $\tilde{f}' \leftarrow \nabla f(\tilde{x})$            ▷ *$n$ queries (non-composition optimization)*
 4:     $x_0 \leftarrow \tilde{x}$
 5:     **for** $k = 0, 1, 2, \ldots, K-1$ **do**
 6:         Uniformly sample pick $i_k$ from $\{1, 2, \ldots, n\}$
 7:         Estimate $\nabla f(x_k)$ using
$$\hat{f}'_k = \tilde{f}' - \nabla F_i(\tilde{x}) + \nabla F_i(x_k)$$
                                                                                   ▷ *2 queries (non-composition optimization)*
 8:         Update $x_{k+1}$ by
$$x_{k+1} \leftarrow x_k - \gamma \hat{f}'_k$$
 9:     **end for**
10:     $\tilde{x}_s \leftarrow x_r$ for randomly chosen $r \in \{0, \cdots, K-1\}$
11: **end for**
---

Note that for strongly convex objectives, the number of inner iterations should be chosen to be in the order of $O(\kappa)$, which leads to a constant linear convergence rate. Thus the query complexity for an error of $\epsilon$ will be

$$O((m + n + m\kappa)\log(1/\epsilon)). \tag{2}$$

## 3 Compositional-SVRG Algorithms

This section introduces two proposed compositional-SVRG algorithms for solving the finite sum composition optimization in (1). In the spirit of SVRG, the two compositional-SVRG algorithms need a reference point $\tilde{x}$ to estimate the gradients (but in different ways). However, unlike SVRG, the estimated gradients are biased due to the "composition" structure in the objective.

### 3.1 Compositional-SVRG-1

In the first proposed algorithm, given a reference point $\tilde{x}$, we first store the gradient $\tilde{f}' = \nabla f(\tilde{x})$ and the value of the inner function $\tilde{G} := G(\tilde{x})$. To estimate the gradient at the current iterate $x_k$, one needs to estimate $G(x_k)$ first by sampling a mini-batch multiset $\mathcal{A}_k$ with size $A$:

$$\hat{G}_k = \tilde{G} - \frac{1}{A} \sum_{1 \leq j \leq A} (G_{\mathcal{A}_k[j]}(\tilde{x}) - G_{\mathcal{A}_k[j]}(x_k)). \tag{3}$$

Based on the estimate of $G(x_k)$, the gradient $\nabla f(x_k)$ is estimated by

$$\hat{f}'_k = (\partial G_{j_k}(x_k))^\top \nabla F_{i_k}(\hat{G}_k) - (\partial G_{j_k}(\tilde{x}))^\top \nabla F_{i_k}(\tilde{G}) + \tilde{f}' \tag{4}$$

where $i_k$ is uniformly sampled from $\{1, 2, \ldots, n\}$ and $j_k$ is uniformly sampled from $\{1, 2, \ldots, m\}$.

Note that unlike SVRG, this estimated $\hat{f}$ is usually *biased*. More specifically,

$$\mathbb{E}_{i_k, j_k, \mathcal{A}_k}(\hat{f}'_k) \neq \nabla f(x_k).$$

This is also the key challenge to prove the linear convergence in the analysis. The Compositional-SVRG-1 algorithm is summarized in Algorithm 2. The query complexity in each step is provided in Algorithm 2 for convenience.

### 3.2 Compositional-SVRG-2

In the second proposed algorithm, given a reference point $\tilde{x}$, we still first store the gradient $\tilde{f}' = \nabla f(\tilde{x})$ and the value of the inner function $\tilde{G} := G(\tilde{x})$. However, here we further store the value of the Jacobian $\tilde{G}' := \partial G(\tilde{x})$. To estimate the gradient at the current iterate $x_k$, one still estimates $G(x_k)$ first by sampling a mini-batch multiset $\mathcal{A}_k$ with size $A$:

$$\hat{G}_k = \tilde{G} - \frac{1}{A} \sum_{1 \leq j \leq A} (G_{\mathcal{A}_k[j]}(\tilde{x}) - G_{\mathcal{A}_k[j]}(x_k)). \tag{5}$$

Here comes the difference from Algorithm 2. We also estimates $\partial G(x_k)$ by sampling a mini-batch multiset $\mathcal{B}_k$ with size $B$:

$$\hat{G}'_k := \tilde{G}' - \frac{1}{B} \sum_{0 \leq j \leq B} \left( \partial G_{\mathcal{B}_k[j]}(\tilde{x}) - \partial G_{\mathcal{B}_k[j]}(x_k) \right). \tag{6}$$

Based on the estimation of $G(x_k)$ and $\partial G(x_k)$, the gradient $\nabla f(x_k)$ is estimated by

$$\hat{f}'_k = (\hat{G}'_k)^\top \nabla F_{i_k}(\hat{G}_k) - (\tilde{G}')^\top \nabla F_{i_k}(\tilde{G}) + \tilde{f}'. \tag{7}$$

where $i_k$ is uniformly sampled from $\{1, 2, \ldots, n\}$. Thus Algorithm 3 features one more estimation in each iteration. This extra computation pays off by an improved convergence rate.

Even though we have an extra estimation here, this estimated $\hat{f}$ is still *biased*. More specifically,

$$\mathbb{E}_{i_k, \mathcal{B}_k, \mathcal{A}_k}(\hat{f}'_k) \neq \nabla f(x_k).$$

The Compositional-SVRG-2 algorithm is summarized in Algorithm 3. The query complexity in each step is provided in Algorithm 3 for convenience.

**Algorithm 2** Compositional-SVRG-1

**Require:** $K$ (the total number of iterations in the inner loop), $S$ (the total number of iterations in the outer loop), $A$ (the size of the minibatch multiset), $\gamma$ (steplength), and $\tilde{x}_0$ (initial point).
**Ensure:** $\tilde{x}_S$.
1: **for** $s = 1, 2, \ldots, S$ **do**
2:     Update the reference point: $\tilde{x} \leftarrow \tilde{x}_{s-1}$
3:     $\tilde{G} \leftarrow G(\tilde{x})$      ▷ $m$ queries
4:     $\tilde{f}' \leftarrow \nabla f(\tilde{x})$      ▷ $m+n$ queries
5:     $x_0 \leftarrow \tilde{x}$
6:     **for** $k = 0, 1, 2, \ldots, K-1$ **do**
7:         Uniformly sample from $\{1, 2, \ldots, m\}$ for $A$ times *with replacement* to form a mini-batch multiset $\mathcal{A}_k$
8:         Estimate $G(x_k)$ by $\hat{G}_k$ using (3)      ▷ $2A$ queries
9:         Uniformly sample $i_k$ from $\{1, 2, \ldots, n\}$ and $j_k$ from $\{1, 2, \ldots, m\}$
10:        Estimate $\nabla f(x_k)$ by $\hat{f}'_k$ using (4)      ▷ 4 queries
11:        Update $x_{k+1}$ by
$$x_{k+1} = x_k - \gamma \hat{f}'_k$$
12:     **end for**
13:     $\tilde{x}_s \leftarrow x_r$ for randomly chosen $r \in \{0, \cdots, K-1\}$
14: **end for**

## 4 Theoretical Analysis

In this section we will show the convergence results for Algorithms 2 and 3. Due to the page limitation, all proofs are provided in the supplement. Before we show the main results, let us make some global assumptions below, which are commonly used for the analysis of stochastic composition optimization algorithms.

**Strongly Convex Objective** $f(x)$ in (1) is strongly convex with parameter $\mu_f$:

$$f(y) \geq f(x) + \langle \nabla f(x), y-x \rangle + \frac{\mu_f}{2}\|y-x\|^2, \quad \forall x, y.$$

**Bounded Jacobian of Inner Functions** We assume that the following upper bounds on all inner component functions:

$$\|\partial G_j(x)\| \leq B_G, \quad \forall x, \forall j \in \{1, \cdots, m\}. \tag{8}$$

**Lipschitzian Gradients** We assume there exist constants $L_F, L_G$ and $L_f$ satisfying $\forall x, \forall y, \forall i \in \{1, \cdots, n\}, \forall j \in \{1, \cdots, m\}$:

$$\|\nabla F_i(x) - \nabla F_i(y)\| \leq L_F \|x-y\|, \tag{9}$$

**Algorithm 3** Compositional-SVRG-2

**Require:** $K, S, A, B, \gamma, \tilde{x}_0$ ▷ The meaning of the variables are the same as in Algorithm 2. $B$ is the size of another minibatch multiset.
**Ensure:** $\tilde{x}_S$.
 1: **for** $s = 1, 2, \ldots, S$ **do**
 2:     Update the reference point $\tilde{x} \leftarrow \tilde{x}_{s-1}$
 3:     $\tilde{G} \leftarrow G(\tilde{x})$     ▷ $m$ queries
 4:     $\tilde{G}' \leftarrow \partial G(\tilde{x})$     ▷ $m$ queries
 5:     $\tilde{f}' \leftarrow \nabla f(\tilde{x})$     ▷ $n$ queries
 6:     $x_0 \leftarrow \tilde{x}$
 7:     **for** $k = 0, 1, 2, \ldots, K-1$ **do**
 8:         Uniformly sample from $\{1, 2, \ldots, m\}$ for $A$ and $B$ times *with replacement* to form two mini-batch multiset $\mathcal{A}_k$ and $\mathcal{B}_k$ respectively
 9:         Estimate $G(x_k)$ by $\hat{G}_k$ using (5)     ▷ $2A$ queries
10:         Estimate $\partial G(x_k)$ by $\hat{G}'_k$ using (6)     ▷ $2B$ queries
11:         Uniformly sample pick $i_k$ from $\{1, 2, \ldots, n\}$
12:         Estimate $\nabla f(x_k)$ by $\hat{f}'_k$ using (7)     ▷ 2 queries
13:         Update $x_{k+1}$ by
$$x_{k+1} \leftarrow x_k - \gamma \hat{f}'_k$$
14:     **end for**
15:     $\tilde{x}_s \leftarrow x_r$ for randomly chosen $r \in \{0, \cdots, K-1\}$
16: **end for**

$$\|\partial G_j(x) - \partial G_j(y)\| \leq L_G \|x - y\|, \tag{10}$$

$$\|(\partial G_j(x))^\top \nabla F_i(G(x)) - (\partial G_j(y))^\top \nabla F_i(G(y))\| \leq L_f \|x - y\|, \tag{11}$$

Note that we immediately have

$$\|\nabla f(x) - \nabla f(y)\| = \frac{1}{mn} \left\| \sum_{i,j} \left( \partial G_j(x) \right)^\top \nabla F_i(G(x)) - (\partial G_j(y))^\top \nabla F_i(G(y)) \right\|$$

$$\leq \frac{1}{mn} \sum_{i,j} \left\| \partial G_j(x))^\top \nabla F_i(G(x)) - (\partial G_j(y))^\top \nabla F_i(G(y)) \right\|$$

$$\leq L_f \|x - y\|, \forall x, y. \tag{12}$$

Thus the Lipschitz constant for the gradient of the whole objective $f$ is also $L_f$.

**Solution Existence** The problem (1) has at least one solution $x^*$.

Next we use two theorems to show our main results on the compositional-SVRG algorithms. All theorems and corollaries hold under the assumptions stated above.

**Theorem 1** (Convergence for Algorithm 2). *For Algorithm 2 we have*

$$\frac{1}{K}\sum_{k=0}^{K-1}\mathbb{E}\|x_k - x^*\|^2 \leqslant \frac{\beta_1}{\beta_2}\mathbb{E}\|\tilde{x} - x^*\|^2,$$

*where*

$$\beta_1 = \frac{1}{K} + \left(\frac{16\gamma B_G^2 L_F^2}{\mu_f} + 4\gamma^2 B_G^2 L_F^2\right)\frac{8B_G^2}{A} + 10\gamma^2 L_f^2;$$

$$\beta_2 = \frac{7\mu_f \gamma}{4} - \left(\frac{16\gamma B_G^2 L_F^2}{\mu_f} + 4\gamma^2 B_G^2 L_F^2\right)\frac{8B_G^2}{A} - 8\gamma^2 L_f^2.$$

To ensure the convergence, one should appropriately choose parameters $\gamma$, $A$, and $K$ to make the ratio $\beta_1/\beta_2 < 1$. The following provides one specification for those parameters.

**Corollary 1** (Linear Rate for Algorithm 2). *Choose parameters in Algorithm 2 as follows:*

$$\gamma = \frac{\mu_f}{32 L_f^2},$$

$$A = \frac{512 B_G^4 L_F^2}{\mu_f^2},$$

$$K = \frac{512 L_f^2}{\mu_f^2},$$

*The following convergence rate for Algorithm 2 holds:*

$$\mathbb{E}\|\tilde{x}_{s+1} - x^*\|^2 = \frac{1}{K}\sum_{k=0}^{K-1}\mathbb{E}\|x_k - x^*\|^2 \leqslant \frac{7}{8}\mathbb{E}\|\tilde{x}_s - x^*\|^2.$$

Corollary 1 essentially suggests a linear convergence rate. To achieve a fixed solution accuracy $\epsilon$, that is, $\mathbb{E}(\|\tilde{x}_s - x^*\|^2) \leq \epsilon$, the required number of queries is $O\left((m + n + KA + K)\log(1/\epsilon)\right) = O\left(\left(m + n + \frac{L_F^2 L_f^2}{\mu_f^4} + \frac{L_f^2}{\mu_f^2}\right)\log(1/\epsilon)\right)$ based on the query complexity of each step in Algorithm 2. Let $\tilde{\kappa}_1 = \max\{\frac{L_F}{\mu_f}, \frac{L_f}{\mu_f}\}$.[3] we see the query complexity of Algorithm 2 is $O((m + n + \tilde{\kappa}_1^4)\log(1/\epsilon))$. Note that this will be smaller than the query complexity (2) for SVRG when $m$ is large. Also note that this $\tilde{\kappa}_1^4$ is much smaller in some special cases, because $L_F/\mu_f$ (and also $L_G/L_f$ in Corollary 2 we will discuss later) could be much smaller than $L_f/\mu_f$.

To analyze the convergence of Algorithm 3, we need one more assumption on the gradients of the outer functions:

---

[3]Note that in classical SVRG, $G$ is an identity function, so $L_F = L_f$. The $\tilde{\kappa}_1$ reduces to the conventional condition number.

**Bounded Gradients of Outer Functions**

$$\|\nabla F_i(x)\| \leq B_F, \forall x, \forall i \in \{1, \cdots, n\}. \tag{13}$$

Then together with the new assumption in (13), we have the following convergence result for Algorithm 3.

**Theorem 2** (Convergence for Algorithm 3). *For algorithm 3 we have*

$$\frac{1}{K} \sum_{k=0}^{K-1} \mathbb{E}(f(x_k) - f^*) \leq \frac{\beta_3}{\beta_4} \mathbb{E}(f(\tilde{x}) - f^*),$$

*where*

$$\beta_3 = \frac{2}{K\mu_f} + \frac{256\gamma B_G^4 L_F^2}{\mu_f A} + \gamma^2 \left( \frac{128}{\mu_f} \left( \frac{B_F^2 L_G^2}{B} + \frac{B_G^4 L_F^2}{A} \right) + 20 L_f \right);$$

$$\beta_4 = \frac{3\gamma}{2} - \frac{256\gamma B_G^4 L_F^2}{\mu_f A} - \gamma^2 \left( \frac{128}{\mu_f} \left( \frac{B_F^2 L_G^2}{B} + \frac{B_G^4 L_F^2}{A} \right) + 16 L_f \right).$$

This theorem admits a similar structure to Theorem 1. We essentially need to appropriately choose $\gamma$, $K$, $A$, and $B$ to make $\beta_3/\beta_4 < 1$. The following corollary provides a specification for these parameters.

**Corollary 2** (Linear Rate for Algorithm 3). *Choose parameters in Algorithm 3 as follows:*

$$\gamma = \frac{1}{320 L_f};$$

$$K \geq \frac{5120 L_f}{\mu_f};$$

$$A \geq \max \left\{ \frac{1024 B_G^4 L_F^2}{\mu_f^2}, \frac{32 B_G^4 L_F^2}{5\mu_f L_f} \right\};$$

$$B \geq \frac{32 B_F^2 L_G^2}{5\mu_f L_f},$$

*we have the following linear convergence rate for Algorithm 3:*

$$\mathbb{E}(f(\tilde{x}_{s+1}) - f^*) = \frac{1}{K} \sum_{k=0}^{K-1} \mathbb{E}(f(x_k) - f^*) \leq \frac{9}{17} \mathbb{E}(f(\tilde{x}_s) - f^*).$$

Let $\tilde{\kappa}_2 = \max\{\frac{L_F}{\mu_f}, \frac{L_G}{L_f}, \frac{L_f}{\mu_f}\}$. Corollary 2 suggests a total query complexity of $O((m + n + K(A + B))\log(1/\epsilon)) = O((m + n + \tilde{\kappa}_2^3)\log(1/\epsilon))$. Note that this will be smaller than the query complexity (2) for SVRG when $m$ is large. Here the query complexity is slightly better[4] than that in Corollary 1. However we need a new assumption that the gradient of $F_i$'s to be bounded and we need an extra estimation for the Jacobian of $G$ at the beginning of each epoch.

---

[4] Here we mean "roughly better", because the definitions of $\kappa$ are different in Corollary 1 and 2, though in most cases they are of the same order.

# 5 Experiment

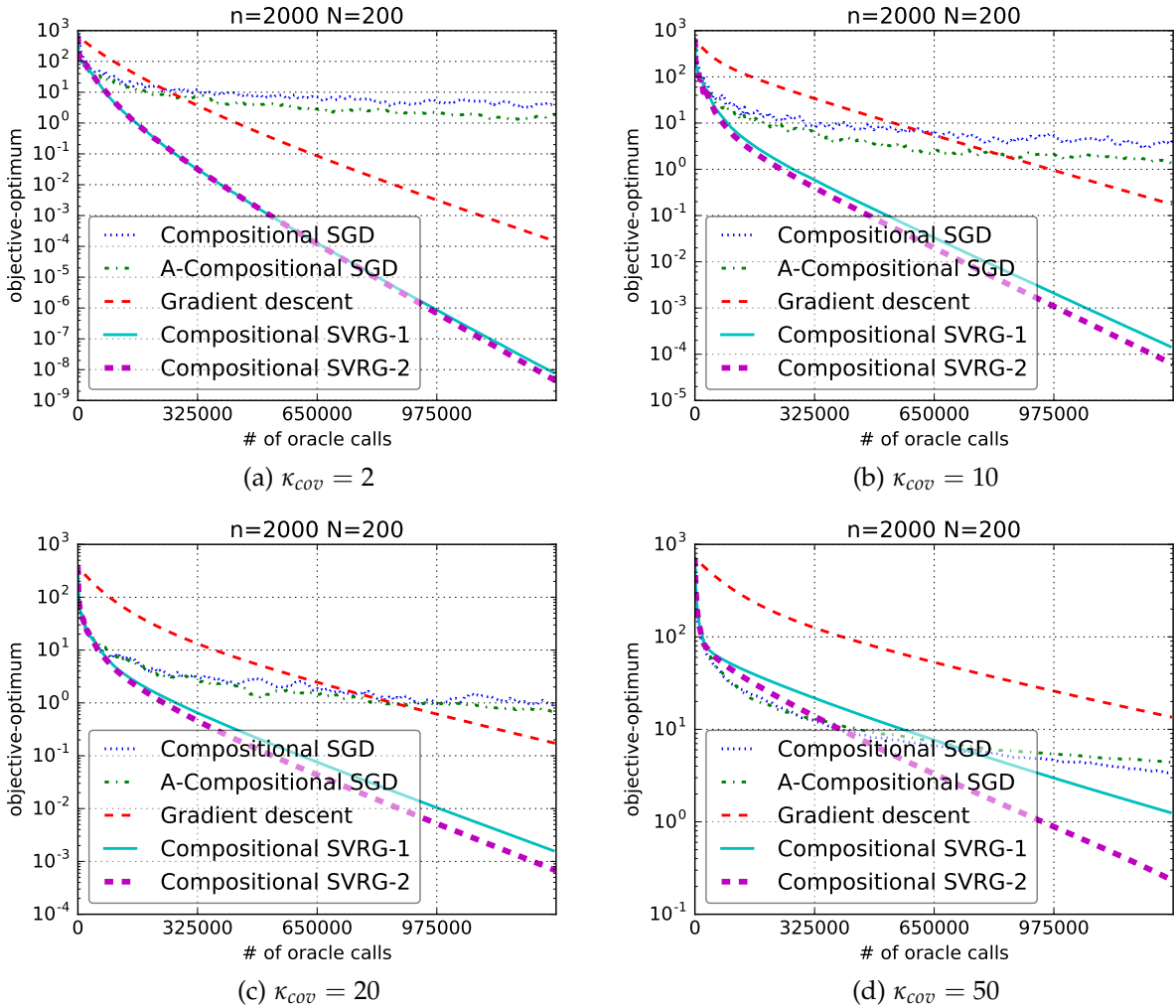

Figure 1: Mean-variance portfolio optimization on synthetic data ($n = 2000, N = 200$). The $y$-axis is the objective value minus the optimal value of the objective. The $x$-axis is the number of oracle calls. The "Compositional SVRG-1" is the Algorithm 2, the "Compositional SVRG-2" is the Algorithm 3. The "Compositional SGD" is the Algorithm 1 in Wang et al. [2014] and The "A-Compositional SGD" is the Algorithm 1 in Wang et al. [2016]. Both SVRG version algorithms use "Compositional-SGD" algorithm to initialize first several steps. The $\kappa_{cov}$ is the conditional number of the covariance matrix of the corresponding Gaussian distribution used to generate reward vectors in each figure. Subfigures (a), (b), (c), and (d) draw the convergence curves for all algorithms with each figure having a different $\kappa_{cov}$.

We conduct empirical studies for the proposed two algorithms by comparing them to three state-of-the-art algorithms:

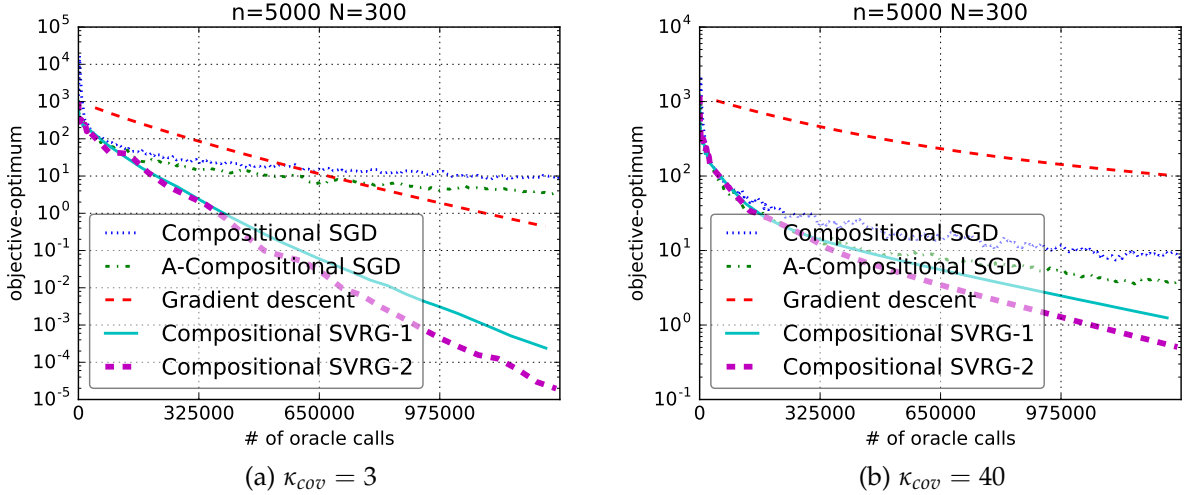

Figure 2: Mean-variance portfolio optimization on synthetic data with $n = 5000, N = 300$. Other settings are the same as Figure 1.

- Gradient descent,
- Compositional SGD [Wang et al., 2014, Algorithm 1],
- Accelerating Compositional SGD [Wang et al., 2016, Algorithm 1].

We use the mean-variance optimization in portfolio management as the objective. Given $N$ assets, let $r_t \in \mathbb{R}^N (t = 1, \ldots, n)$ be the reward vectors observed at different time points. The goal is to maximize the return of the investment as well as controlling the investment risk. Let $x \in \mathbb{R}^N$ be the quantities invested to each portfolio. The problem can be formulated into the following mean-variance optimization[5]:

$$\max_x \quad \frac{1}{n}\sum_{i=1}^n \langle r_i, x \rangle - \frac{1}{n}\sum_{i=1}^n \left( \langle r_i, x \rangle - \frac{1}{n}\sum_{j=1}^n \langle r_j, x \rangle \right)^2.$$

It can be viewed as an instance of the composition optimization (1) with the specification for $G_j(\cdot)$ and $F_i(\cdot)$ as the following:

$$G_j(x) = \begin{pmatrix} x \\ \langle r_j, x \rangle \end{pmatrix}, j = 1, \ldots, n;$$

$$F_i(y) = -y_{N+1} + (\langle r_i, y_{1:N} \rangle - y_{N+1})^2, i = 1, \ldots, n,$$

where $y_{1:N}$ denotes the sub-vector consisting of the first $N$ components of the vector $y \in \mathbb{R}^{N+1}$ and $y_{N+1}$ denotes the last component of $y$.

---

[5]This formulation is just used for proof of concept. A more efficient way would be simply calculating $\sum_j r_j$. Then the problem reduces to a standard stochastic optimization.

In the experiment we choose $n = 2000$ and $N = 200$. The reward vectors are generated with the procedure below:

1. Generate a random Gaussian distribution on $\mathbb{R}^N$ with the condition number of its covariance matrix denoted by $\kappa_{cov}$.

2. Each $r_i$ is sampled from this Gaussian distribution with all elements set to its absolute value to make sure the problem has a solution.

We can then control $\kappa_{cov}$ to roughly control the $\tilde{\kappa}_1$ and $\tilde{\kappa}_2$ of our composition optimization problem, because $\tilde{\kappa}_1$ and $\tilde{\kappa}_2$ are proportional to $\kappa_{cov}$. We report the comparison results in Figure 1. The initial points are chosen to be the same for all algorithms and the x-axis in Figure 1 is the computational cost measured by the number of queries to the oracle. That is, whenever the algorithm queries $\nabla F_i(y)$ or $\partial G_i(x)$ or $G_i(x)$ for some $i$ at some point, the x-axis value is incremented by 1. Like the SVRG algorithm [Johnson and Zhang, 2013], both compositional-SVRG algorithms run the compositional-SGD algorithm (which is the Algorithm 1 in Wang et al. [2014]) for the first 10000 iterations and then run the proposed SVRG algorithms. We observe that

- The proposed two algorithms (Compositional SVRG-1 and Compositional SVRG-2) converge at a linear rate and outperform other algorithms overall;

- Compositional SVRG-2 becomes faster than Compositional SVRG-1 when $\kappa_{cov}$ becomes larger, while they are comparable when $\kappa_{cov}$ are small. This observation is consistent with our theoretical analysis, since Compositional SVRG-2 has a better dependency on the condition number than Compositional SVRG-1.

We also test our algorithms on problem with a larger size ($n = 5000, N = 300$), and show the results in Figure 2.

# 6 Conclusion and Future Work

This paper considers the finite-sum composition optimization and proposes two efficient algorithm by using the SVRG technique to reduce variance of compositional gradient. The proposed two algorithms admit the linear convergence rate for strongly convex objectives with query complexity $O((m+n+\tilde{\kappa}_1^4)\log(1/\epsilon))$ and $O((m+n+\tilde{\kappa}_2^3)\log(1/\epsilon))$ respectively. To the best of our knowledge, this is the first work to study the general finite-sum composition optimization. The future work will be 1) the convergence rate and query complexity for weakly convex problem; 2) the convergence rate and query complexity for nonconvex optimization; 3) how (or is it possible) to improve the query complexity to $O((m+n+\kappa)\log(1/\epsilon))$ to make it consistent with SVRG for the classical stochastic optimization?

**Acknowledgements**

Mengdi Wang was supported by NSF CMMI 1653435 and NSF DMS 1619818. We thank the reviewers for their constructive comments and suggestions, especially for pointing out the potentially more efficient implementation in the experiment.

# SUPPLEMENTARY MATERIAL

**Proof to Theorem 1**

*Proof.* We start from decomposing the expectation of $\|x_{k+1} - x^*\|^2$:

$$\begin{aligned}
&\mathbb{E}\|x_{k+1} - x^*\|^2 \\
&= \mathbb{E}\|x_k - x^*\|^2 + \mathbb{E}\|x_{k+1} - x_k\|^2 + 2\mathbb{E}\langle x_{k+1} - x_k, x_k - x^*\rangle \\
&= \mathbb{E}\|x_k - x^*\|^2 + \gamma^2 \mathbb{E}\left\|(\partial G_{j_k}(x_k))^\top \nabla F_{i_k}(\hat{G}_k) - (\partial G_{j_k}(\tilde{x}))^\top \nabla F_{i_k}(\tilde{G}) + \nabla f(\tilde{x})\right\|^2 \\
&\quad - 2\gamma \mathbb{E}\left\langle (\partial G_{j_k}(x_k))^\top \nabla F_{i_k}(\hat{G}_k) - (\partial G_{j_k}(\tilde{x}))^\top \nabla F_{i_k}(\tilde{G}) + \nabla f(\tilde{x}), x_k - x^*\right\rangle. 
\end{aligned} \quad (14)$$

Given the observation

$$\begin{aligned}
&\mathbb{E}\left\langle -(\partial G_{j_k}(\tilde{x}))^\top \nabla F_{i_k}(\tilde{G}) + \nabla f(\tilde{x}), x_k - x^*\right\rangle \\
&= \mathbb{E}\left\langle -\mathbb{E}_{j_k,i_k} J_{G_{j_k}}^\top(\tilde{x})\nabla F_{i_k}(\tilde{G}) + \nabla f(\tilde{x}), x_k - x^*\right\rangle \\
&= \mathbb{E}\langle -\nabla f(\tilde{x}) + \nabla f(\tilde{x}), x_k - x^*\rangle \\
&= 0,
\end{aligned}$$

It follows from (14) that

$$\begin{aligned}
\mathbb{E}\|x_{k+1} - x^*\|^2 &= \mathbb{E}\|x_k - x^*\|^2 - 2\gamma \underbrace{\mathbb{E}\left\langle (\partial G_{j_k}(x_k))^\top \nabla F_{i_k}(\hat{G}_k), x_k - x^*\right\rangle}_{=:T_1} \\
&\quad + \gamma^2 \underbrace{\mathbb{E}\left\|(\partial G_{j_k}(x_k))^\top \nabla F_{i_k}(\hat{G}_k) - (\partial G_{j_k}(\tilde{x}))^\top \nabla F_{i_k}(\tilde{G}) + \nabla f(\tilde{x})\right\|^2}_{=:T_2}.
\end{aligned} \quad (15)$$

We then bound $T_1$. From the strong convexity of $f(x)$ we have the following inequality:

$$\langle \nabla f(x), x - x^*\rangle \geq \mu_f \|x_k - x^*\|^2. \quad (16)$$

It follows that

$$\begin{aligned}
T_1 &= \mathbb{E}\left\langle (\partial G_{j_k}(x_k))^\top \nabla F_{i_k}(\hat{G}_k), x_k - x^*\right\rangle \\
&= \mathbb{E}\left\langle (\partial G_{j_k}(x_k))^\top \nabla F_{i_k}(\hat{G}_k) - \nabla f(x_k), x_k - x^*\right\rangle + \mathbb{E}\langle \nabla f(x_k), x_k - x^*\rangle \\
&\stackrel{(16)}{\geq} \underbrace{\mathbb{E}\left\langle (\partial G_{j_k}(x_k))^\top \nabla F_{i_k}(\hat{G}_k) - \nabla f(x_k), x_k - x^*\right\rangle}_{=:T_3} + \mathbb{E}\mu_f\|x_k - x^*\|^2. 
\end{aligned} \quad (17)$$

We then bound $T_3$. Recall that for any $\alpha > 0$ we have

$$\frac{1}{\alpha}x^2 + \alpha y^2 \geq 2|\langle x, y\rangle| \geq |\langle x, y\rangle|. \quad (18)$$

It follows that

$$
\begin{aligned}
T_3 &= \mathbb{E}\left\langle (\partial G_{j_k}(x_k))^\top \nabla F_{i_k}(\hat{G}_k) - \nabla f(x_k), x_k - x^* \right\rangle \\
&= \mathbb{E}\left\langle (\partial G_{j_k}(x_k))^\top \nabla F_{i_k}(\hat{G}_k) - (\partial G_{j_k}(x_k))^\top \nabla F_{i_k}(G(x_k)), x_k - x^* \right\rangle \\
&\overset{(18)}{\geq} -\frac{1}{\alpha}\underbrace{\mathbb{E}\left\| (\partial G_{j_k}(x_k))^\top \nabla F_{i_k}(\hat{G}_k) - (\partial G_{j_k}(x_k))^\top \nabla F_{i_k}(G(x_k)) \right\|^2}_{=:T_4} - \alpha \mathbb{E}\|x_k - x^*\|^2, \forall \alpha > 0.
\end{aligned}
\tag{19}
$$

For $T_4$, from the definition of $\hat{G}_k$,

$$
\begin{aligned}
T_4 &= \mathbb{E}\left\| (\partial G_{j_k}(x_k))^\top \nabla F_{i_k}(\hat{G}_k) - (\partial G_{j_k}(x_k))^\top \nabla F_{i_k}(G(x_k)) \right\|^2 \\
&= \mathbb{E}\left\| (\partial G_{j_k}(x_k))^\top \nabla F_{i_k}\left( \tilde{G} - \frac{1}{A}\sum_{1 \leq j \leq A}(G_{\mathcal{A}_k[j]}(\tilde{x}) - G_{\mathcal{A}_k[j]}(x_k)) \right) - (\partial G_{j_k}(x_k))^\top \nabla F_{i_k}(G(x_k)) \right\|^2 \\
&\leq \mathbb{E}\left\| (\partial G_{j_k}(x_k))^\top \right\|^2 \left\| \nabla F_{i_k}\left( \tilde{G} - \frac{1}{A}\sum_{1 \leq j \leq A}(G_{\mathcal{A}_k[j]}(\tilde{x}) - G_{\mathcal{A}_k[j]}(x_k)) \right) - \nabla F_{i_k}(G(x_k)) \right\|^2 \\
&\overset{(8)}{\leq} B_G^2 \mathbb{E}\left\| \nabla F_{i_k}\left( \tilde{G} - \frac{1}{A}\sum_{1 \leq j \leq A}(G_{\mathcal{A}_k[j]}(\tilde{x}) - G_{\mathcal{A}_k[j]}(x_k)) \right) - \nabla F_{i_k}(G(x_k)) \right\|^2 \\
&\overset{(9)}{\leq} B_G^2 L_F^2 \underbrace{\mathbb{E}\left\| \tilde{G} - \frac{1}{A}\sum_{1 \leq j \leq A}(G_{\mathcal{A}_k[j]}(\tilde{x}) - G_{\mathcal{A}_k[j]}(x_k)) - G(x_k) \right\|^2}_{=:T_0}.
\end{aligned}
\tag{20}
$$

Let $\alpha = \frac{\mu_f}{8}$ in (19) and put the bound of $T_4$ in it, we obtain

$$
\begin{aligned}
T_3 &\overset{(19)}{\geq} -\frac{1}{\alpha}T_4 - \alpha \mathbb{E}\|x_k - x^*\|^2 \\
&\overset{(20)}{\geq} -\frac{8B_G^2 L_F^2}{\mu_f}T_0 - \frac{\mu_f}{8}\mathbb{E}\|x_k - x^*\|^2.
\end{aligned}
\tag{21}
$$

Then put this bound on $T_3$ to (17).

$$
\begin{aligned}
T_1 &\overset{(17)}{\geq} T_3 + \mathbb{E}\mu_f \|x_k - x^*\|^2 \\
&\overset{(21)}{\geq} -\frac{8B_G^2 L_F^2}{\mu_f}T_0 + \frac{7\mu_f}{8}\mathbb{E}\|x_k - x^*\|^2.
\end{aligned}
\tag{22}
$$

Now we have $T_1$ bounded. We use this bound to bound the $T_1$ in the equality (15) at the beginning.

$$
\begin{aligned}
\mathbb{E}\|x_{k+1} - x^*\|^2 &\overset{(15)}{=} \mathbb{E}\|x_k - x^*\|^2 - 2\gamma T_1 + \gamma^2 T_2 \\
&\overset{(22)}{\leq} \mathbb{E}\|x_k - x^*\|^2 - \frac{7\mu_f \gamma}{4}\mathbb{E}\|x_k - x^*\|^2 + \frac{16\gamma B_G^2 L_F^2}{\mu_f}T_0 + \gamma^2 T_2.
\end{aligned}
\tag{23}
$$

We then bound $T_2$. Recall for any $\beta$ we have

$$\|\beta_1 + \beta_2 + \cdots + \beta_t\|^2 \leqslant t \left(\|\beta_1\|^2 + \cdots + \|\beta_t\|^2\right), \forall t \in \mathbb{N}_+. \tag{24}$$

From the definition of $T_2$ in (15) we have the following bound on $T_2$:

$$\begin{aligned}
T_2 &= \mathbb{E}\left\|(\partial G_{j_k}(x_k))^\top \nabla F_{i_k}(\hat{G}_k) - (\partial G_{j_k}(\tilde{x}))^\top \nabla F_{i_k}(\tilde{G}) + \nabla f(\tilde{x})\right\|^2 \\
&\stackrel{(24)}{\leqslant} 2\mathbb{E}\|\nabla f(\tilde{x})\|^2 + 2\mathbb{E}\left\|(\partial G_{j_k}(x_k))^\top \nabla F_{i_k}(\hat{G}_k) - (\partial G_{j_k}(\tilde{x}))^\top \nabla F_{i_k}(\tilde{G})\right\|^2 \\
&\stackrel{(24)}{\leqslant} 2\mathbb{E}\|\nabla f(\tilde{x})\|^2 + 4\underbrace{\mathbb{E}\left\|(\partial G_{j_k}(x_k))^\top \nabla F_{i_k}(\hat{G}_k) - (\partial G_{j_k}(x_k))^\top \nabla F_{i_k}(G(x_k))\right\|^2}_{\text{the same as } T_4} \\
&\quad + 4\underbrace{\mathbb{E}\left\|(\partial G_{j_k}(x_k))^\top \nabla F_{i_k}(G(x_k)) - (\partial G_{j_k}(\tilde{x}))^\top \nabla F_{i_k}(\tilde{G})\right\|^2}_{=:T_5} \\
&\stackrel{(20)}{\leqslant} 2\mathbb{E}\|\nabla f(\tilde{x})\|^2 + 4B_G^2 L_F^2 T_0 + 4T_5. \tag{25}
\end{aligned}$$

To bound $T_5$, we simply use the Lipschitzian condition (11)

$$\begin{aligned}
T_5 &= \mathbb{E}\left\|(\partial G_{j_k}(x_k))^\top \nabla F_{i_k}(g(x_k)) - (\partial G_{j_k}(\tilde{x}))^\top \nabla F_{i_k}(\tilde{G})\right\|^2 \\
&\stackrel{(11)}{\leqslant} L_f^2 \mathbb{E}\|x_k - \tilde{x}\|^2,
\end{aligned}$$

Put this bound back to (25) we obtain

$$T_2 \stackrel{(25)}{\leqslant} 2\mathbb{E}\|\nabla f(\tilde{x})\|^2 + 4B_G^2 L_F^2 T_0 + 4L_f^2 \mathbb{E}\|x_k - \tilde{x}\|^2. \tag{26}$$

Now we have $T_2$ bounded, and we put this bound back to (23).

$$\begin{aligned}
\mathbb{E}\|x_{k+1} - x^*\|^2 &\stackrel{(23)}{\leqslant} \mathbb{E}\|x_k - x^*\|^2 - \frac{7\mu_f \gamma}{4}\mathbb{E}\|x_k - x^*\|^2 + \frac{16\gamma B_G^2 L_F^2}{\mu_f} T_0 + \gamma^2 T_2 \\
&\stackrel{(26)}{\leqslant} \mathbb{E}\|x_k - x^*\|^2 - \frac{7\mu_f \gamma}{4}\mathbb{E}\|x_k - x^*\|^2 + \frac{16\gamma B_G^2 L_F^2}{\mu_f} T_0 \\
&\quad + 2\gamma^2 \mathbb{E}\|\nabla f(\tilde{x})\|^2 + 4\gamma^2 B_G^2 L_F^2 T_0 + 4\gamma^2 L_f^2 \mathbb{E}\|x_k - \tilde{x}\|^2 \\
&\stackrel{(12)}{=} \mathbb{E}\|x_k - x^*\|^2 - \frac{7\mu_f \gamma}{4}\mathbb{E}\|x_k - x^*\|^2 + 2\gamma^2 L_f^2 \mathbb{E}\|\tilde{x} - x^*\|^2 \\
&\quad + \left(\frac{16\gamma B_G^2 L_F^2}{\mu_f} + 4\gamma^2 B_G^2 L_F^2\right) T_0 + 4\gamma^2 L_f^2 \mathbb{E}\|x_k - \tilde{x}\|^2, \tag{27}
\end{aligned}$$

where the last step comes from (12) by letting $x = x_k$ and $y = x^*$.

There is still one term, $T_0$, not bounded. We now start to bound it. From the definition of $T_0$ in (20):

$$T_0 = \mathbb{E}\left\|\tilde{G} - \frac{1}{A}\sum_{1 \leq j \leq A}(G_{\mathcal{A}_k[j]}(\tilde{x}) - G_{\mathcal{A}_k[j]}(x_k)) - G(x_k)\right\|^2$$

$$
\begin{aligned}
&= \mathbb{E}\left\|\frac{1}{A}\sum_{1\leq j\leq A}(G_{\mathcal{A}_k[j]}(\tilde{x}) - G_{\mathcal{A}_k[j]}(x_k)) - (\tilde{G} - G(x_k))\right\|^2 \\
&= \frac{1}{A^2}\mathbb{E}\left\|\sum_{1\leq j\leq A}((G_{\mathcal{A}_k[j]}(\tilde{x}) - G_{\mathcal{A}_k[j]}(x_k)) - (\tilde{G} - G(x_k)))\right\|^2 \\
&= \frac{1}{A^2}\sum_{1\leq j\leq A}\mathbb{E}\|(G_{\mathcal{A}_k[j]}(\tilde{x}) - G_{\mathcal{A}_k[j]}(x_k)) - (\tilde{G} - G(x_k))\|^2,
\end{aligned}
$$

where the last step comes from the fact that the indices in $\mathcal{A}_k$ are independent. Specifically,

$$
\begin{aligned}
&\mathbb{E}\left\|\sum_{1\leq j\leq A}(G_{\mathcal{A}_k[j]}(\tilde{x}) - G_{\mathcal{A}_k[j]}(x_k) - \tilde{G} + G(x_k))\right\|^2 \\
&= \mathbb{E}\sum_{1\leq j\leq A}\|(G_{\mathcal{A}_k[j]}(\tilde{x}) - G_{\mathcal{A}_k[j]}(x_k) - \tilde{G} + G(x_k))\|^2 \\
&\quad + 2\mathbb{E}\sum_{1\leq j'<j\leq A}\langle(G_{\mathcal{A}_k[j]}(\tilde{x}) - G_{\mathcal{A}_k[j]}(x_k) - \tilde{G} + G(x_k)), (G_{\mathcal{A}_k[j']}(\tilde{x}) - G_{\mathcal{A}_k[j']}(x_k) - \tilde{G} + G(x_k))\rangle \\
&= \mathbb{E}\sum_{1\leq j\leq A}\|(G_{\mathcal{A}_k[j]}(\tilde{x}) - G_{\mathcal{A}_k[j]}(x_k) - \tilde{G} + G(x_k))\|^2 \\
&\quad + 2\mathbb{E}\sum_{1\leq j'<j\leq A}\langle\mathbb{E}_{\mathcal{A}_k[j]}(G_{\mathcal{A}_k[j]}(\tilde{x}) - G_{\mathcal{A}_k[j]}(x_k) - \tilde{G} + G(x_k)), (G_{\mathcal{A}_k[j']}(\tilde{x}) - G_{\mathcal{A}_k[j']}(x_k) - \tilde{G} + G(x_k))\rangle \\
&= \mathbb{E}\sum_{1\leq j\leq A}\|(G_{\mathcal{A}_k[j]}(\tilde{x}) - G_{\mathcal{A}_k[j]}(x_k) - \tilde{G} + G(x_k))\|^2 \\
&\quad + 2\mathbb{E}\sum_{1\leq j'<j\leq A}\langle 0, (G_{\mathcal{A}_k[j']}(\tilde{x}) - G_{\mathcal{A}_k[j']}(x_k) - \tilde{G} + G(x_k))\rangle \\
&= \mathbb{E}\sum_{1\leq j\leq A}\|(G_{\mathcal{A}_k[j]}(\tilde{x}) - G_{\mathcal{A}_k[j]}(x_k) - \tilde{G} + G(x_k))\|^2.
\end{aligned}
\tag{28}
$$

Finally $T_0$ can be bounded by

$$
\begin{aligned}
T_0 &= \frac{1}{A^2}\sum_{1\leq j\leq A}\mathbb{E}\|(G_{\mathcal{A}_k[j]}(\tilde{x}) - G_{\mathcal{A}_k[j]}(x_k)) - (\tilde{G} - G(x_k))\|^2 \\
&\stackrel{(24)}{\leq} \frac{4}{A^2}\sum_{1\leq j\leq A}\mathbb{E}\left(\|G_{\mathcal{A}_k[j]}(\tilde{x}) - G_{\mathcal{A}_k[j]}(x^*)\|^2 + \left\|G_{\mathcal{A}_k[j]}(x_k) - G_{\mathcal{A}_k[j]}(x^*)\right\|^2 \right. \\
&\qquad\qquad \left. + \|\tilde{G} - G(x^*)\|^2 + \|G(x_k) - G(x^*)\|^2\right) \\
&\stackrel{(8)}{\leq} \frac{8B_G^2}{A^2}\sum_{1\leq j\leq A}\mathbb{E}(\|\tilde{x} - x^*\|^2 + \|x_k - x^*\|^2) \\
&= \frac{8B_G^2}{A}\mathbb{E}(\|\tilde{x} - x^*\|^2 + \|x_k - x^*\|^2).
\end{aligned}
\tag{29}
$$

By passing this bound to (27) we finally get all $T$ terms bounded:

$$\mathbb{E}\|x_{k+1} - x^*\|^2$$

$$\stackrel{(27)}{\leqslant} \mathbb{E}\|x_k - x^*\|^2 - \frac{7\mu_f\gamma}{4}\mathbb{E}\|x_k - x^*\|^2 + 2\gamma^2 L_f^2 \mathbb{E}\|\tilde{x} - x^*\|^2$$
$$+ \left(\frac{16\gamma B_G^2 L_F^2}{\mu_f} + 4\gamma^2 B_G^2 L_F^2\right) T_0 + 4\gamma^2 L_f^2 \mathbb{E}\|x_k - \tilde{x}\|^2$$

$$\stackrel{(29)}{\leqslant} \mathbb{E}\|x_k - x^*\|^2 - \frac{7\mu_f\gamma}{4}\mathbb{E}\|x_k - x^*\|^2 + 2\gamma^2 L_f^2 \mathbb{E}\|\tilde{x} - x^*\|^2$$
$$+ \left(\frac{16\gamma B_G^2 L_F^2}{\mu_f} + 4\gamma^2 B_G^2 L_F^2\right) \frac{8B_G^2}{A}\mathbb{E}(\|\tilde{x} - x^*\|^2 + \|x_k - x^*\|^2) + 4\gamma^2 L_f^2 \mathbb{E}\|x_k - x^* + x^* - \tilde{x}\|^2$$

$$\stackrel{(24)}{\leqslant} \mathbb{E}\|x_k - x^*\|^2 - \frac{7\mu_f\gamma}{4}\mathbb{E}\|x_k - x^*\|^2 + 2\gamma^2 L_f^2 \mathbb{E}\|\tilde{x} - x^*\|^2$$
$$+ \left(\frac{16\gamma B_G^2 L_F^2}{\mu_f} + 4\gamma^2 B_G^2 L_F^2\right) \frac{8B_G^2}{A}\mathbb{E}(\|\tilde{x} - x^*\|^2 + \|x_k - x^*\|^2)$$
$$+ 8\gamma^2 L_f^2 \mathbb{E}(\|x_k - x^*\|^2 + \|\tilde{x} - x^*\|^2)$$
$$= \mathbb{E}\|x_k - x^*\|^2$$
$$- \left(\frac{7\mu_f\gamma}{4} - \left(\frac{16\gamma B_G^2 L_F^2}{\mu_f} + 4\gamma^2 B_G^2 L_F^2\right) \frac{8B_G^2}{A} - 8\gamma^2 L_f^2\right) \mathbb{E}\|x_k - x^*\|^2$$
$$+ \left(\left(\frac{16\gamma B_G^2 L_F^2}{\mu_f} + 4\gamma^2 B_G^2 L_F^2\right) \frac{8B_G^2}{A} + 10\gamma^2 L_f^2\right) \mathbb{E}\|\tilde{x} - x^*\|^2.$$

Summing this inequality from $k = 0$ to $k = K - 1$, we obtain

$$\mathbb{E}\|x_{k+K} - x^*\|^2$$
$$\leqslant \mathbb{E}\|\tilde{x} - x^*\|^2$$
$$- \left(\frac{7\mu_f\gamma}{4} - \left(\frac{16\gamma B_G^2 L_F^2}{\mu_f} + 4\gamma^2 B_G^2 L_F^2\right) \frac{8B_G^2}{A} - 8\gamma^2 L_f^2\right) \sum_{k=0}^{K-1}\mathbb{E}\|x_k - x^*\|^2$$
$$+ K\left(\left(\frac{16\gamma B_G^2 L_F^2}{\mu_f} + 4\gamma^2 B_G^2 L_F^2\right) \frac{8B_G^2}{A} + 10\gamma^2 L_f^2\right) \mathbb{E}\|\tilde{x} - x^*\|^2.$$

Discarding the left hand side, we complete the proof by

$$\frac{1}{K}\sum_{k=0}^{K-1}\mathbb{E}\|x_k - x^*\|^2$$
$$\leqslant \frac{\frac{1}{K} + \left(\frac{16\gamma B_G^2 L_F^2}{\mu_f} + 4\gamma^2 B_G^2 L_F^2\right) \frac{8B_G^2}{A} + 10\gamma^2 L_f^2}{\frac{7\mu_f\gamma}{4} - \left(\frac{16\gamma B_G^2 L_F^2}{\mu_f} + 4\gamma^2 B_G^2 L_F^2\right) \frac{8B_G^2}{A} - 8\gamma^2 L_f^2}\mathbb{E}\|\tilde{x} - x^*\|^2.$$

$\square$

**Proof to Corollary 1**

*Proof.* To appropriately choose $\gamma, K$ and $A$ in Algorithm 2, the key is to ensure the coefficient $\frac{\beta_1}{\beta_2} < 1$ in Theorem 1:

$$\frac{\beta_1}{\beta_2} = \frac{\frac{1}{K} + \left(\frac{16\gamma B_G^2 L_F^2}{\mu_f} + 4\gamma^2 B_G^2 L_F^2\right)\frac{8B_G^2}{A} + 10\gamma^2 L_f^2}{\frac{7\mu_f \gamma}{4} - \left(\frac{16\gamma B_G^2 L_F^2}{\mu_f} + 4\gamma^2 B_G^2 L_F^2\right)\frac{8B_G^2}{A} - 8\gamma^2 L_f^2}.$$

We choose $A$ satisfying both

$$4\gamma^2 B_G^2 L_F^2 \frac{8B_G^2}{A} \leq \frac{\mu_f \gamma}{4},$$

$$\frac{16\gamma B_G^2 L_F^2}{\mu_f} \frac{8B_G^2}{A} \leq \frac{\mu_f \gamma}{4},$$

which is equivalent to

$$A \geq \max\left\{\frac{128\gamma B_G^4 L_F^2}{\mu_f}, \frac{512 B_G^4 L_F^2}{\mu_f^2}\right\}.$$

We choose $\gamma$ satisfying

$$8\gamma^2 L_f^2 \leq \frac{\mu_f \gamma}{4},$$

which is equivalent to

$$\gamma \leq \frac{\mu_f}{32 L_f^2}.$$

It follows that

$$\frac{\frac{1}{K} + \left(\frac{16\gamma B_G^2 L_F^2}{\mu_f} + 4\gamma^2 B_G^2 L_F^2\right)\frac{8B_G^2}{A} + 10\gamma^2 L_f^2}{\frac{7\mu_f \gamma}{4} - \left(\frac{16\gamma B_G^2 L_F^2}{\mu_f} + 4\gamma^2 B_G^2 L_F^2\right)\frac{8B_G^2}{A} - 8\gamma^2 L_f^2}$$

$$\leq \frac{\frac{1}{K} + \frac{13\mu_f \gamma}{16}}{\mu_f \gamma}$$

$$= \frac{13}{16} + \frac{1}{K\mu_f \gamma}.$$

We then choose $K$ satisfying

$$\frac{1}{K\mu_f \gamma} \leq \frac{1}{16},$$

which is equivalent to

$$K \geq \frac{16}{\mu_f \gamma}.$$

Thus choosing $\gamma$, $A$, and $K$ appropriately in the following to satisfy all conditions derived above

$$\gamma = \frac{\mu_f}{32L_f^2},$$

$$A = \frac{512 B_G^4 L_F^2}{\mu_f^2},$$

$$K = \frac{512 L_f^2}{\mu_f^2},$$

we obtain a linear convergence rate of coefficient $\frac{\beta_1}{\beta_2} = \frac{7}{8}$ from Theorem 1. □

**Lemma 1.** *Under the assumption in* (11), *we have*

$$\frac{1}{n} \sum_{i=1}^{n} \left\| \frac{1}{m} \sum_{j=1}^{m} (\partial G_j(x))^\top \nabla F_i(G(x)) - \frac{1}{m} \sum_{j=1}^{m} (\partial G_j(x^*))^\top \nabla F_i(G(x^*)) \right\|^2 \leq 2L_f(f(x) - f^*).$$

*Proof.* Recall that at the optimal point we always have

$$f'(x^*) = \frac{1}{mn} \sum_{i=1}^{n} \sum_{j=1}^{m} (\partial G_j(x^*))^\top \nabla F_i(G(x^*)) = 0. \tag{30}$$

We can derive the Lipschitz constant of $F_i(G(x))$ from (11)

$$\|\nabla F_i(G(x)) - \nabla F_i(G(y))\|$$
$$= \frac{1}{m} \left\| \sum_j (\partial G_j(x))^\top \nabla F_i(G(x)) - \sum_j (\partial G_j(y))^\top \nabla F_i(G(y)) \right\|$$
$$\leq \frac{1}{m} \sum_j \|(\partial G_j(x))^\top \nabla F_i(G(x)) - (\partial G_j(y))^\top \nabla F_i(G(y))\|$$
$$\leq L_f \|x - y\|, \forall i. \tag{31}$$

From this Lipschitz condition, we obtain

$$F_i(G(x)) \overset{(31)}{\geq} F_i(G(x^*)) + \frac{1}{m} \left\langle \sum_{j=1}^{m} (\partial G_j(x^*))^\top \nabla F_i(G(x^*)), x - x^* \right\rangle$$
$$+ \frac{1}{2L_f} \left\| \frac{1}{m} \sum_{j=1}^{m} (\partial G_j(x))^\top \nabla F_i(x) - \frac{1}{m} \sum_{j=1}^{m} (\partial G_j(x^*))^\top \nabla F_i(x^*) \right\|^2.$$

Summing from $i = 1$ to $i = n$, using (30) and noting that $\frac{1}{n} \sum_{i=1}^{n} F_i(G(x)) = f(x)$, we obtain

$$\frac{1}{n} \sum_{i=1}^{n} \left\| \frac{1}{m} \sum_{j=1}^{m} (\partial G_j(x))^\top \nabla F_i(x) - \frac{1}{m} \sum_{j=1}^{m} (\partial G_j(x^*))^\top \nabla F_i(x^*) \right\|^2 \leq 2L_f(f(x) - f^*),$$

completing the proof. □

**Proof to Theorem 2**

*Proof.* Note that in this proof we redefine the terms $T_1, T_2, \ldots$, and they may not refer to the same expressions in the proof of Theorem 1. From

$$x_{k+1} - x_k = -\gamma((\hat{G}'_k)^\top \nabla F_{i_k}(\hat{G}_k) - (\tilde{G}')^\top \nabla F_{i_k}(\tilde{G}) + \tilde{f}').$$

we immediately obtain

$$\begin{aligned} \mathbb{E}\|x_{k+1} - x^*\|^2 &= \mathbb{E}\|x_k - x^*\|^2 + \mathbb{E}\|x_{k+1} - x_k\|^2 + 2\mathbb{E}\langle x_{k+1} - x_k, x_k - x^*\rangle \\ &= \mathbb{E}\|x_k - x^*\|^2 + \gamma^2 \mathbb{E}\|(\hat{G}'_k)^\top \nabla F_{i_k}(\hat{G}_k) - (\tilde{G}')^\top \nabla F_{i_k}(\tilde{G}) + \tilde{f}'\|^2 \\ &\quad - 2\gamma \mathbb{E}\langle (\hat{G}'_k)^\top \nabla F_{i_k}(\hat{G}_k) - (\tilde{G}')^\top \nabla F_{i_k}(\tilde{G}) + \tilde{f}', x_k - x^*\rangle. \end{aligned}$$

Note that the last term can be simplified:

$$\begin{aligned} &\mathbb{E}\langle (\hat{G}'_k)^\top \nabla F_{i_k}(\hat{G}_k) - (\tilde{G}')^\top \nabla F_{i_k}(\tilde{G}) + \tilde{f}', x_k - x^*\rangle \\ &= \mathbb{E}\langle (\hat{G}'_k)^\top \nabla F_{i_k}(\hat{G}_k) - \mathbb{E}_{i_k}(\tilde{G}')^\top \nabla F_{i_k}(\tilde{G}) + \tilde{f}', x_k - x^*\rangle \\ &= \mathbb{E}\langle (\hat{G}'_k)^\top \nabla F_{i_k}(\hat{G}_k) - \tilde{f}' + \tilde{f}', x_k - x^*\rangle \\ &= \mathbb{E}\langle (\hat{G}'_k)^\top \nabla F_{i_k}(\hat{G}_k), x_k - x^*\rangle. \end{aligned}$$

Therefore, we have

$$\begin{aligned} \mathbb{E}\|x_{k+1} - x^*\|^2 &= \mathbb{E}\|x_k - x^*\|^2 - 2\gamma \underbrace{\mathbb{E}\langle (\hat{G}'_k)^\top \nabla F_{i_k}(\hat{G}_k), x_k - x^*\rangle}_{=:T_1} \\ &\quad + \gamma^2 \underbrace{\mathbb{E}\|(\hat{G}'_k)^\top \nabla F_{i_k}(\hat{G}_k) - (\tilde{G}')^\top \nabla F_{i_k}(\tilde{G}) + \tilde{f}'\|^2}_{=:T_2}. \end{aligned} \quad (32)$$

First we estimate the lower bound for $T_1$:

$$\begin{aligned} T_1 &= \mathbb{E}\langle (\hat{G}'_k)^\top \nabla F_{i_k}(\hat{G}_k), x_k - x^*\rangle \\ &= \underbrace{\mathbb{E}\langle (\hat{G}'_k)^\top \nabla F_{i_k}(\hat{G}_k) - \nabla f(x_k), x_k - x^*\rangle}_{=:T_3} + \mathbb{E}\langle \nabla f(x_k), x_k - x^*\rangle \\ &\geqslant T_3 + \mathbb{E}(f(x_k) - f^*). \end{aligned} \quad (33)$$

Then we estimate the lower bound for $T_3$

$$\begin{aligned} T_3 &= \mathbb{E}\langle (\hat{G}'_k)^\top \nabla F_{i_k}(\hat{G}_k) - \nabla f(x_k), x_k - x^*\rangle \\ &= \mathbb{E}\left\langle (\partial G_{j_k}(x_k))^\top \nabla F_{i_k}(\hat{G}_k) - (\partial G_{j_k}(x_k))^\top \nabla F_{i_k}(G(x_k)), x_k - x^*\right\rangle, \end{aligned}$$

where $j_k$ is a new (imaginary) random variable that is chosen uniformly randomly from $\{1, \cdots, m\}$ and is independent of other random variables. $\mathbb{E}$ also takes expectation on $j_k$. Thus using the same technique as we use in (19) while proving Theorem 1, we obtain

$$T_3 \geqslant -\frac{1}{\alpha} \underbrace{\mathbb{E}\|(\partial G_{j_k}(x_k))^\top \nabla F_{i_k}(\hat{G}_k) - (\partial G_{j_k}(x_k))^\top \nabla F_{i_k}(G(x_k))\|^2}_{=:T_4} - \alpha \mathbb{E}\|x_k - x^*\|^2, \forall \alpha > 0.$$

and

$$T_4 = \mathbb{E}\|(\partial G_{j_k}(x_k))^\top \nabla F_{i_k}(\hat{G}_k) - (\partial G_{j_k}(x_k))^\top \nabla F_{i_k}(G(x_k))\|^2$$

$$\stackrel{(20)}{\leqslant} B_G^2 L_F^2 T_0,$$

where

$$T_0 := \mathbb{E}\left\|\tilde{G} - \frac{1}{A}\sum_{1\leq j\leq A}(G_{\mathcal{A}_k[j]}(\tilde{x}) - G_{\mathcal{A}_k[j]}(x_k)) - G(x_k)\right\|^2.$$

Let $\alpha = \frac{\mu_f}{8}$, we obtain

$$T_3 \geqslant -\frac{8B_G^2 L_F^2}{\mu_f}T_0 - \frac{\mu_f}{8}\mathbb{E}\|x_k - x^*\|^2.$$

Put the bound of $T_3$ into (33) and note that

$$\mu_f \|x_k - x^*\|^2 \leqslant 2(f(x_k) - f^*). \tag{34}$$

We obtain

$$T_1 \geqslant -\frac{8B_G^2 L_F^2}{\mu_f}T_0 - \frac{\mu_f}{8}\mathbb{E}\|x_k - x^*\|^2 + \mathbb{E}(f(x_k) - f^*)$$

$$\stackrel{(34)}{\geqslant} -\frac{8B_G^2 L_F^2}{\mu_f}T_0 + \frac{3}{4}\mathbb{E}(f(x_k) - f^*). \tag{35}$$

Now we have $T_1$ bounded. We then start to bound $T_2$. From the definition of $T_2$ we have

$$T_2 = \mathbb{E}\|(\hat{G}_k')^\top \nabla F_{i_k}(\hat{G}_k) - (\tilde{G}')^\top \nabla F_{i_k}(\tilde{G}) + \tilde{f}'\|^2$$

$$\stackrel{(24)}{\leqslant} 2\mathbb{E}\|\tilde{f}'\|^2 + 2\mathbb{E}\|(\hat{G}_k')^\top \nabla F_{i_k}(\hat{G}_k) - (\tilde{G}')^\top \nabla F_{i_k}(\tilde{G})\|^2$$

$$= 2\mathbb{E}\|\tilde{f}'\|^2 + 2\mathbb{E}\left\|(\hat{G}_k')^\top \nabla F_{i_k}(\hat{G}_k) - \frac{1}{m}\sum_{j=1}^m (\partial G_j(x_k))^\top \nabla F_{i_k}(G(x_k))\right. \tag{36}$$

$$\left. + \frac{1}{m}\sum_{j=1}^m (\partial G_j(x_k))^\top \nabla F_{i_k}(G(x_k)) - (\tilde{G}')^\top \nabla F_{i_k}(\tilde{G})\right\|^2$$

$$\stackrel{(24)}{\leqslant} 2\mathbb{E}\|\tilde{f}'\|^2 + 4\mathbb{E}\left\|(\hat{G}_k')^\top \nabla F_{i_k}(\hat{G}_k) - \frac{1}{m}\sum_{j=1}^m (\partial G_j(x_k))^\top \nabla F_{i_k}(G(x_k))\right\|^2$$

$$+ 4\mathbb{E}\left\|\frac{1}{m}\sum_{j=1}^m (\partial G_j(x_k))^\top \nabla F_{i_k}(G(x_k)) - (\tilde{G}')^\top \nabla F_{i_k}(\tilde{G})\right\|^2$$

$$\stackrel{(24)}{\leqslant} 2\mathbb{E}\|\tilde{f}'\|^2 + 4\mathbb{E}\left\|(\hat{G}_k')^\top \nabla F_{i_k}(\hat{G}_k) - \frac{1}{m}\sum_{j=1}^m (\partial G_j(x_k))^\top \nabla F_{i_k}(G(x_k))\right\|^2$$

$$+ 8\mathbb{E}\left\|\frac{1}{m}\sum_{j=1}^m (\partial G_j(x_k))^\top \nabla F_{i_k}(G(x_k)) - \frac{1}{m}\sum_{j=1}^m (\partial G_j(x^*))^\top \nabla F_{i_k}(G(x^*))\right\|^2$$

$$+8\mathbb{E}\left\|(\tilde{G}')^\top \nabla F_{i_k}(\tilde{G}) - \frac{1}{m}\sum_{j=1}^m (\partial G_j(x^*))^\top \nabla F_{i_k}(G(x^*))\right\|^2$$

$$\leq 4L_f \mathbb{E}(f(\tilde{x}) - f^*) + 4\underbrace{\mathbb{E}\left\|(\hat{G}'_k)^\top \nabla F_{i_k}(\hat{G}_k) - \frac{1}{m}\sum_{j=1}^m (\partial G_j(x_k))^\top \nabla F_{i_k}(G(x_k))\right\|^2}_{=:T_5}$$

$$+16(L_f \mathbb{E}(f(\tilde{x}) - f^*) + L_f \mathbb{E}(f(x_k) - f^*)), \tag{37}$$

where the last step comes from Lemma 1 and $\frac{\|\nabla f(\tilde{x})\|^2}{2L_f} \leq f(\tilde{x}) - f^*$.

Note that $T_5$ can be bounded by

$$T_5 = \mathbb{E}\left\|(\hat{G}'_k)^\top \nabla F_{i_k}(\hat{G}_k) - \frac{1}{m}\sum_{j=1}^m (\partial G_j(x_k))^\top \nabla F_{i_k}(G(x_k))\right\|^2$$

$$\overset{(24)}{\leq} 2\mathbb{E}\left\|(\hat{G}'_k)^\top \nabla F_{i_k}(\hat{G}_k) - \frac{1}{m}\sum_{j=1}^m (\partial G_j(x_k))^\top \nabla F_{i_k}(\hat{G}_k)\right\|^2$$

$$+2\mathbb{E}\left\|\frac{1}{m}\sum_{j=1}^m (\partial G_j(x_k))^\top \nabla F_{i_k}(G(x_k)) - \frac{1}{m}\sum_{j=1}^m (\partial G_j(x_k))^\top \nabla F_{i_k}(\hat{G}_k)\right\|^2$$

$$\overset{(13),(8)}{\leq} 2B_F^2 \mathbb{E}\left\|(\hat{G}'_k)^\top - \frac{1}{m}\sum_{j=1}^m (\partial G_j(x_k))^\top\right\|^2 + 2B_G^2 \mathbb{E}\|\nabla F_{i_k}(G(x_k)) - \nabla F_{i_k}(\hat{G}_k)\|^2$$

$$\leq 2B_F^2 \mathbb{E}\left\|(\tilde{G}')^\top - \frac{1}{B}\left(\sum_{1\leq j\leq B}((\partial G_{\mathcal{B}_j[j]}(\tilde{x}))^\top - (\partial G_{\mathcal{B}_k[j]}(x_k))^\top)\right) - \frac{1}{m}\sum_{j=1}^m (\partial G_j(x_k))^\top\right\|^2$$

$$+2B_G^2 L_F^2 \mathbb{E}\left\|\tilde{G} - \frac{1}{A}\sum_{1\leq j\leq A}(G_{\mathcal{A}_k[j]}(\tilde{x}) - G_{\mathcal{A}_k[j]}(x_k)) - G(x_k)\right\|^2$$

$$= \frac{2B_F^2}{B^2}\mathbb{E}\left\|-\sum_{1\leq j\leq B}\left(((\partial G_{\mathcal{B}_k[j]}(\tilde{x}))^\top - (\partial G_{\mathcal{B}_k[j]}(x_k))^\top) - \left((\tilde{G}')^\top - \frac{1}{m}\sum_{j=1}^m (\partial G_j(x_k))^\top\right)\right)\right\|^2$$

$$+\frac{2B_G^2 L_F^2}{A^2}\mathbb{E}\left\|-\sum_{1\leq j\leq A}((G_{\mathcal{A}_k[j]}(\tilde{x}) - G_{\mathcal{A}_k[j]}(x_k)) - (G(x_k) - \tilde{G}))\right\|^2$$

Using the same technique as in (28), the above inequality continues as

$$= \frac{2B_F^2}{B^2}\mathbb{E}\sum_{1\leq j\leq B}\left\|(\partial G_{\mathcal{B}_k[j]}(x_k))^\top - (\partial G_{\mathcal{B}_k[j]}(\tilde{x}))^\top + \left((\tilde{G}')^\top - \frac{1}{m}\sum_{j=1}^m (\partial G_j(x_k))^\top\right)\right\|^2$$

$$+\frac{2B_G^2 L_F^2}{A^2}\mathbb{E}\sum_{1\leq j\leq A}\|G_{\mathcal{A}_k[j]}(x_k) - G_{\mathcal{A}_k[j]}(\tilde{x}) + (G(x_k) - \tilde{G})\|^2$$

$$\overset{(8),(10),(24)}{\leq} \frac{2B_F^2}{B^2}\sum_{1\leq j\leq B} 8L_G^2(\mathbb{E}\|\tilde{x} - x^*\|^2 + \mathbb{E}\|x_k - x^*\|^2)$$

$$+ \frac{2B_G^2 L_F^2}{A^2} \sum_{1 \leq j \leq A} 8B_G^2 (\mathbb{E}\|\tilde{x} - x^*\|^2 + \mathbb{E}\|x_k - x^*\|^2)$$

$$= 16 \left( \frac{B_F^2 L_G^2}{B} + \frac{B_G^4 L_F^2}{A} \right) (\mathbb{E}\|\tilde{x} - x^*\|^2 + \mathbb{E}\|x_k - x^*\|^2)$$

$$\leq \frac{32}{\mu_f} \left( \frac{B_F^2 L_G^2}{B} + \frac{B_G^4 L_F^2}{A} \right) (\mathbb{E}(f(\tilde{x}) - f^*) + \mathbb{E}(f(x_k) - f^*)).$$

Now we continue to bound $T_2$ in (37) using the bound for $T_5$ above:

$$T_2 \stackrel{(37)}{\leq} 4L_f \mathbb{E}(f(\tilde{x}) - f^*) + 4T_5 + 16(L_f(f(\tilde{x}) - f^*) + L_f \mathbb{E}(f(x_k) - f^*))$$
$$= 20L_f \mathbb{E}(f(\tilde{x}) - f^*) + 16L_f \mathbb{E}(f(x_k) - f^*) + 4T_5$$
$$\leq 20L_f \mathbb{E}(f(\tilde{x}) - f^*) + 16L_f \mathbb{E}(f(x_k) - f^*)$$
$$+ \frac{128}{\mu_f} \left( \frac{B_F^2 L_G^2}{B} + \frac{B_G^4 L_F^2}{A} \right) (\mathbb{E}(f(\tilde{x}) - f^*) + \mathbb{E}(f(x_k) - f^*))$$
$$= \left( \frac{128}{\mu_f} \left( \frac{B_F^2 L_G^2}{B} + \frac{B_G^4 L_F^2}{A} \right) + 16L_f \right) \mathbb{E}(f(x_k) - f^*)$$
$$+ \left( \frac{128}{\mu_f} \left( \frac{B_F^2 L_G^2}{B} + \frac{B_G^4 L_F^2}{A} \right) + 20L_f \right) \mathbb{E}(f(\tilde{x}) - f^*). \tag{38}$$

Now we have $T_2$ bounded. Finally we put the bounds of $T_2, T_1$ in (38) and (35) into (32) and note that using the same procedure in the proof of Theorem 1 (see (29)) we have

$$T_0 \leq \frac{8B_G^2}{A} \mathbb{E}(\|\tilde{x} - x^*\|^2 + \|x_k - x^*\|^2). \tag{39}$$

We obtain:

$$\mathbb{E}\|x_{k+1} - x^*\|^2$$
$$\stackrel{(32)}{=} \mathbb{E}\|x_k - x^*\|^2 - 2\gamma T_1 + \gamma^2 T_2$$
$$\stackrel{(35),(38)}{\leq} \mathbb{E}\|x_k - x^*\|^2 - 2\gamma \left( -\frac{8B_G^2 L_F^2}{\mu_f} T_0 + \frac{3}{4} \mathbb{E}(f(x_k) - f^*) \right)$$
$$+ \gamma^2 \left( \frac{128}{\mu_f} \left( \frac{B_F^2 L_G^2}{B} + \frac{B_G^4 L_F^2}{A} \right) + 16L_f \right) \mathbb{E}(f(x_k) - f^*)$$
$$+ \gamma^2 \left( \frac{128}{\mu_f} \left( \frac{B_F^2 L_G^2}{B} + \frac{B_G^4 L_F^2}{A} \right) + 20L_f \right) \mathbb{E}(f(\tilde{x}) - f^*)$$
$$\stackrel{(39)}{\leq} \mathbb{E}\|x_k - x^*\|^2 - 2\gamma \left( -\frac{8B_G^2 L_F^2}{\mu_f} \left( \frac{8B_G^2}{A} \mathbb{E}(\|\tilde{x} - x^*\|^2 + \|x_k - x^*\|^2) \right) + \frac{3}{4} \mathbb{E}(f(x_k) - f^*) \right)$$
$$+ \gamma^2 \left( \frac{128}{\mu_f} \left( \frac{B_F^2 L_G^2}{B} + \frac{B_G^4 L_F^2}{A} \right) + 16L_f \right) \mathbb{E}(f(x_k) - f^*)$$

$$
\begin{aligned}
&+\gamma^2 \left( \frac{128}{\mu_f} \left( \frac{B_F^2 L_G^2}{B} + \frac{B_G^4 L_F^2}{A} \right) + 20 L_f \right) \mathbb{E}(f(\tilde{x}) - f^*) \\
\overset{(34)}{\leqslant}\;& \mathbb{E}\|x_k - x^*\|^2 - 2\gamma \left( -\frac{128 B_G^4 L_F^2}{\mu_f^2 A} \mathbb{E}(f(\tilde{x}) - f^* + f(x_k) - f^*) + \frac{3}{4} \mathbb{E}(f(x_k) - f^*) \right) \\
&+\gamma^2 \left( \frac{128}{\mu_f} \left( \frac{B_F^2 L_G^2}{B} + \frac{B_G^4 L_F^2}{A} \right) + 16 L_f \right) \mathbb{E}(f(x_k) - f^*) \\
&+\gamma^2 \left( \frac{128}{\mu_f} \left( \frac{B_F^2 L_G^2}{B} + \frac{B_G^4 L_F^2}{A} \right) + 20 L_f \right) \mathbb{E}(f(\tilde{x}) - f^*) \\
=\;& \mathbb{E}\|x_k - x^*\|^2 \\
&- \left( \frac{3\gamma}{2} - \frac{256 \gamma B_G^4 L_F^2}{\mu_f^2 A} - \gamma^2 \left( \frac{128}{\mu_f} \left( \frac{B_F^2 L_G^2}{B} + \frac{B_G^4 L_F^2}{A} \right) + 16 L_f \right) \right) \mathbb{E}(f(x_k) - f^*) \\
&+ \left( \frac{256 \gamma B_G^4 L_F^2}{\mu_f^2 A} + \gamma^2 \left( \frac{128}{\mu_f} \left( \frac{B_F^2 L_G^2}{B} + \frac{B_G^4 L_F^2}{A} \right) + 20 L_f \right) \right) \mathbb{E}(f(\tilde{x}) - f^*).
\end{aligned}
$$

Summing from $k = 0$ to $k = K - 1$, we obtain

$$
\begin{aligned}
\mathbb{E}\|x_K - x^*\|^2 \leqslant\;& \mathbb{E}\|\tilde{x} - x^*\|^2 \\
&- \left( \frac{3\gamma}{2} - \frac{256 \gamma B_G^4 L_F^2}{\mu_f^2 A} - \gamma^2 \left( \frac{128}{\mu_f} \left( \frac{B_F^2 L_G^2}{B} + \frac{B_G^4 L_F^2}{A} \right) + 16 L_f \right) \right) \sum_{k=0}^{K-1} \mathbb{E}(f(x_k) - f^*) \\
&+ \left( \frac{256 \gamma B_G^4 L_F^2}{\mu_f^2 A} + \gamma^2 \left( \frac{128}{\mu_f} \left( \frac{B_F^2 L_G^2}{B} + \frac{B_G^4 L_F^2}{A} \right) + 20 L_f \right) \right) K \mathbb{E}(f(\tilde{x}) - f^*).
\end{aligned}
$$

Discarding the LHS and note that $\|\tilde{x} - x^*\|^2 \leqslant \frac{2}{\mu_f}(f(\tilde{x}) - f^*)$, we obtain

$$
\frac{1}{K} \sum_{k=0}^{K-1} \mathbb{E}(f(x_k) - f^*) \leqslant \frac{\frac{2}{K\mu_f} + \frac{256 \gamma B_G^4 L_F^2}{\mu_f^2 A} + \gamma^2 \left( \frac{128}{\mu_f} \left( \frac{B_F^2 L_G^2}{B} + \frac{B_G^4 L_F^2}{A} \right) + 20 L_f \right)}{\frac{3\gamma}{2} - \frac{256 \gamma B_G^4 L_F^2}{\mu_f^2 A} - \gamma^2 \left( \frac{128}{\mu_f} \left( \frac{B_F^2 L_G^2}{B} + \frac{B_G^4 L_F^2}{A} \right) + 16 L_f \right)} \mathbb{E}(f(\tilde{x}) - f^*),
$$

completing the proof. □

**Proof to Corollary 2**

*Proof.* To appropriately choose parameters $\gamma$, $K$, $A$, and $B$, the key is to ensure the coefficient $\frac{\beta_3}{\beta_4} < 1$ in Therom 2:

$$
\frac{\beta_3}{\beta_4} = \frac{\frac{2}{K\mu_f} + \frac{256 \gamma B_G^4 L_F^2}{\mu_f^2 A} + \gamma^2 \left( \frac{128}{\mu_f} \left( \frac{B_F^2 L_G^2}{B} + \frac{B_G^4 L_F^2}{A} \right) + 20 L_f \right)}{\frac{3\gamma}{2} - \frac{256 \gamma B_G^4 L_F^2}{\mu_f^2 A} - \gamma^2 \left( \frac{128}{\mu_f} \left( \frac{B_F^2 L_G^2}{B} + \frac{B_G^4 L_F^2}{A} \right) + 16 L_f \right)}.
$$

We choose $A$, $B$, and $\gamma$ satisfying (40), (41), (42), and (43):

$$\frac{256\gamma B_G^4 L_F^2}{\mu_f^2 A} \leq \frac{\gamma}{4} \tag{40}$$

$$\Rightarrow A \geq \frac{1024 B_G^4 L_F^2}{\mu_f^2}$$

$$\gamma^2 \frac{128}{\mu_f} \frac{B_F^2 L_G^2}{B} \leq \frac{\gamma}{16} \tag{41}$$

$$\Rightarrow B \geq \gamma \frac{2048}{\mu_f} B_F^2 L_G^2$$

$$\gamma^2 \frac{128}{\mu_f} \frac{B_G^4 L_F^2}{A} \leq \frac{\gamma}{16} \tag{42}$$

$$\Rightarrow A \geq \gamma \frac{2048}{\mu_f} B_G^4 L_F^2$$

$$20\gamma^2 L_f \leq \frac{\gamma}{16} \tag{43}$$

$$\Rightarrow \gamma \leq \frac{1}{320 L_f}.$$

Then we have the following bound on the coefficient

$$\nu = \frac{\frac{2}{K\mu_f} + \frac{256\gamma B_G^4 L_F^2}{\mu_f A} + \gamma^2 \left( \frac{128}{\mu_f} \left( \frac{B_F^2 L_G^2}{B} + \frac{B_G^4 L_F^2}{A} \right) + 20 L_f \right)}{\frac{3\gamma}{2} - \frac{256\gamma B_G^4 L_F^2}{\mu_f A} - \gamma^2 \left( \frac{128}{\mu_f} \left( \frac{B_F^2 L_G^2}{B} + \frac{B_G^4 L_F^2}{A} \right) + 16 L_f \right)}$$

$$\leq \frac{\frac{2}{K\mu_f} + \frac{\gamma}{4} + \frac{3\gamma}{16}}{\frac{3\gamma}{2} - \frac{\gamma}{4} - \frac{3\gamma}{16}}$$

$$= \frac{\frac{2}{K\mu_f} + \frac{7\gamma}{16}}{\frac{17\gamma}{16}}$$

$$= \frac{32}{17 K \mu_f \gamma} + \frac{7}{17}.$$

We then choose $K$ satisfying

$$\frac{32}{17 K \mu_f \gamma} \leq \frac{2}{17},$$

which is equivalent to

$$K \geq \frac{16}{\mu_f \gamma}.$$

Thus choosing $\gamma$, $A$, and $K$ appropriately in the following to satisfy all conditions derived above

$$\gamma = \frac{1}{320 L_f}$$

$$K \geq \frac{16}{\mu_f \gamma} = \frac{5120 L_f}{\mu_f}$$

$$A \geq \max\left\{ \frac{1024 B_G^4 L_F^2}{\mu_f^2}, \gamma \frac{2048}{\mu_f} B_G^4 L_F^2 \right\} = \max\left\{ \frac{1024 B_G^4 L_F^2}{\mu_f^2}, \frac{32 B_G^4 L_F^2}{5 \mu_f L_f} \right\}$$

$$B \geq \gamma \frac{2048}{\mu_f} B_F^2 L_G^2 = \frac{32 B_F^2 L_G^2}{5 \mu_f L_f},$$

we will obtain a 9/17 linear convergence rate with $\frac{\beta_3}{\beta_4} = \frac{9}{17}$ from Theorem 2. □